\documentclass[11pt]{article}

\usepackage[T1]{fontenc}
\usepackage[margin=1in]{geometry}
\usepackage{amsmath,amssymb,amsthm,microtype}
\usepackage[unicode,hidelinks]{hyperref}
\hypersetup{
  pdftitle={A 19/5 lower bound for the four-dimensional Erdos box problem},
  pdfauthor={Dean Menezes}
}

\newtheorem{theorem}{Theorem}
\newtheorem{lemma}[theorem]{Lemma}
\newtheorem{corollary}[theorem]{Corollary}
\numberwithin{equation}{section}

\DeclareMathOperator{\Tr}{Tr}
\DeclareMathOperator{\N}{N}
\DeclareMathOperator{\Span}{span}
\newcommand{\F}{\mathbf{F}}

\title{A \texorpdfstring{$19/5$}{19/5} lower bound for the
four-dimensional Erd\H{o}s box problem}
\author{Dean Menezes}
\date{}

\begin{document}
\maketitle

\begin{abstract}
For every power $q\ge2$ of~$2$, we construct a box-free subset of
$[q^5]^4$ with $q^{19}$ points, improving the lower-bound exponent
in the four-dimensional Erd\H{o}s box problem from $15/4$ to $19/5$.
In fact, $[q^5]^4$ partitions into $q$ isomorphic box-free sets.
The proof uses a finite-field rectangle lemma and orthogonality in a
four-dimensional bilinear space. We also deduce a uniform bound
$f_4(n)\ge2^{-19/20}n^{19/5}$ and a lower bound for grid-coloring numbers.
\end{abstract}

\section{Introduction}\label{sec:intro}

A \emph{box} in a product $V_1\times\cdots\times V_r$ is a set
$A_1\times\cdots\times A_r$ with $|A_i|=2$ for all~$i$.
A set is \emph{box-free} if it contains no box. Let $f_r(n)$ be the
maximum size of a box-free subset of $[n]^r$, where
$[n]=\{1,\ldots,n\}$. Equivalently, $f_r(n)$ is the maximum number
of edges in an $r$-partite $r$-uniform hypergraph with $n$ vertices
in each part and no copy of $K_{2,\ldots,2}^{(r)}$.

Erd\H{o}s~\cite{Erd64} proved that
$f_r(n)=O(n^{r-1/2^{r-1}})$ and conjectured that this is the correct
order of magnitude. The conjecture is known for $r=3$~\cite{Men26+}.
For $r=4$, Gunderson, R\"odl, and Sidorenko~\cite{GRS99} obtained
the lower bound $f_4(n)=\Omega(n^{15/4})$, subsequently recovered by
Conlon, Pohoata, and Zakharov~\cite{CPZ21} and Gordeev~\cite{Gor24}.
We improve the lower-bound exponent to~$19/5$.

\begin{theorem}\label{thm:main}
For every power $q\ge2$ of~$2$, there is a box-free subset of
$[q^5]^4$ of size~$q^{19}$. Moreover, $[q^5]^4$ partitions into
$q$ isomorphic box-free subsets.
\end{theorem}

The construction uses the inverse of the permutation
$Q(x)=x^{q^2+q}$ of $K=\F_{q^3}$. A rectangle lemma shows that
four points of the form $Q^{-1}(u_i+v_j)$ are either collinear or
affinely independent over $\F_q$. This dichotomy, together with a
bilinear pairing, rules out boxes. We give the lemma and the proof,
then deduce the uniform lower bound and the coloring consequence.

\section{A finite-field rectangle lemma}\label{sec:lemma}

Fix a power $q\ge2$ of~$2$, and write $k=\F_q$ and $K=\F_{q^3}$.
The relative trace and norm are
\[
  \Tr(x)=x+x^q+x^{q^2},
  \qquad
  \N(x)=x^{q^2+q+1}.
\]
Since $\Tr(\alpha)=\alpha$ for $\alpha\in k$, the trace pairing
$(x,y)\mapsto\Tr(xy)$ is nondegenerate: for $x\ne0$, take $y=x^{-1}$.
Also, $\N(x)=0$ if and only if $x=0$.

Set $Q(x)=x^{q^2+q}$. Expanding $\N(x+y)$ and using
$xQ(x)=\N(x)\in k$ gives
\begin{equation}\label{eq:norm_id}
  \Tr\bigl((x+y)(Q(x)+Q(y))\bigr)=\N(x+y)
  \qquad(x,y\in K).
\end{equation}
Thus $Q(x)=Q(y)$ implies $x=y$, so $Q$ permutes~$K$.
Its polarization
\[
  P(x,y)=Q(x+y)+Q(x)+Q(y)=x^qy^{q^2}+x^{q^2}y^q
\]
is symmetric and $k$-bilinear, with $P(x,x)=0$. Applying
\eqref{eq:norm_id} to $x+y$ and $y$ shows that
$\Tr(xP(x,y))=0$. Moreover, $P(x,y)=0$ if and only if $x,y$ are
linearly dependent over~$k$: when $x\ne0$, this equation is equivalent
to $(y/x)^q=y/x$.

The next lemma gives the dichotomy used in the construction. We
include the proof, following~\cite[Lemma~3.2]{Men26+}.

\begin{lemma}\label{lem:rectangle}
Let $u_0,u_1,v_0,v_1\in K$, and put
\[
  A_{ij}=Q^{-1}(u_i+v_j)
  \qquad(i,j\in\{0,1\}).
\]
If $u_0\ne u_1$ and $v_0\ne v_1$, then either the four points $A_{ij}$
are affinely independent over~$k$, or they are collinear and satisfy
\[
  \sum_{i,j}A_{ij}=0,
  \qquad
  \sum_{i,j}\Tr(u_iA_{ij})\ne0.
\]
If $u_0=u_1$ or $v_0=v_1$, both sums vanish.
\end{lemma}

\begin{proof}
If $u_0=u_1$ or $v_0=v_1$, the entries occur in equal pairs and both
sums vanish. Assume that $u_0\ne u_1$ and $v_0\ne v_1$, and set
\[
  a=A_{10}+A_{00},\qquad
  b=A_{01}+A_{00},\qquad
  c=\sum_{i,j}A_{ij}.
\]
Then $a,b\ne0$ and $A_{11}=A_{00}+a+b+c$. Summing
$Q(A_{ij})=u_i+v_j$ and expanding gives
\begin{equation}\label{eq:Q_rect}
  Q(c)+P(A_{00}+a+b,c)+P(a,b)=0.
\end{equation}
Multiplying by $c$ and taking the trace yields
\[
  \N(c)=\Tr\bigl(cP(a,b)\bigr).
\]
The $k$-trilinear form $(a,b,c)\mapsto\Tr(cP(a,b))$ is alternating,
since it vanishes whenever two arguments coincide. If $c\ne0$, then
$\N(c)\ne0$, so $a,b,c$ are linearly independent and the four points
$A_{ij}$ are affinely independent.

If $c=0$, equation~\eqref{eq:Q_rect} gives $P(a,b)=0$. Hence
$b=\lambda a$ for some $\lambda\in k^\times$, and the points are
collinear. Finally, \eqref{eq:norm_id} gives
\begin{align*}
  \sum_{i,j}\Tr(u_iA_{ij})
  &=\Tr\bigl((u_0+u_1)b\bigr)\\
  &=\lambda\Tr\bigl(a(Q(A_{00}+a)+Q(A_{00}))\bigr)\\
  &=\lambda\N(a)\ne0.\qedhere
\end{align*}
\end{proof}

\section{The construction}\label{sec:construction}

Continue with $k,K,Q$ from Section~\ref{sec:lemma}. Equip
$W=K\oplus k$ with the nondegenerate symmetric bilinear form
\[
  \langle(a,\alpha),(b,\beta)\rangle=\Tr(ab)+\alpha\beta,
\]
and put $e=(0,1)$, so $\langle e,e\rangle=1$.
For $x=(u,w)$ and $y=(v,z)$ in~$W$, define the pair-map
$\Phi\colon W\times W\to W$ by
\[
  \Phi(x,y)=\bigl(a,\Tr(au)+w+z\bigr),
  \qquad a=Q^{-1}(u+v).
\]
Let $V=W\times k$, so $|V|=q^5$. Write a vertex as
$v_p=(x_p,\omega_p)$, where $x_p=(u_p,w_p)$. Define
$\chi\colon V^4\to k$ by
\begin{equation}\label{eq:coloring}
\begin{split}
  \chi(v_1,v_2,v_3,v_4)
  ={}&\langle\Phi(x_1,x_2),\Phi(x_3,x_4)\rangle\\
     &+w_1w_2+w_3w_4+\sum_{p=1}^4\omega_p.
\end{split}
\end{equation}

The pairing will make the two difference spaces of a box orthogonal.
The asymmetric trace term in $\Phi$ forces a two-dimensional difference
space to contain~$e$ when both $K$-coordinates vary; the terms
$w_1w_2+w_3w_4$ rule out the case when neither varies.

\begin{proof}[Proof of Theorem~\ref{thm:main}]
In \eqref{eq:coloring}, the first three vertices and $x_4$ determine
$\omega_4$ uniquely for any prescribed value of~$\chi$. Thus every fiber has $q^{19}$ points.
The translation $\omega_4\mapsto\omega_4+\alpha$ carries
$S=\chi^{-1}(0)$ onto $\chi^{-1}(\alpha)$ and preserves boxes. It
therefore suffices to prove that $S$ is box-free.

Suppose that $S$ contains a box with vertices
$v_p^{(i)}=(x_p^{(i)},\omega_p^{(i)})$, where
$x_p^{(i)}=(u_p^{(i)},w_p^{(i)})$, $p\in\{1,2,3,4\}$, and
$i\in\{0,1\}$. Comparing two corners differing only in part~$p$
shows that $x_p^{(0)}=x_p^{(1)}$ would force
$\omega_p^{(0)}=\omega_p^{(1)}$. Since the vertices in each part are
distinct, we have
\begin{equation}\label{eq:distinct}
  x_p^{(0)}\ne x_p^{(1)}
  \qquad(p=1,2,3,4).
\end{equation}

For $i,j,r,s\in\{0,1\}$, put
\[
  C_{ij}=\Phi(x_1^{(i)},x_2^{(j)}),
  \qquad C'_{rs}=\Phi(x_3^{(r)},x_4^{(s)}),
\]
and define
\begin{align*}
  L&=\Span_k\{C_{ij}+C_{00}:i,j\in\{0,1\}\},\\
  L'&=\Span_k\{C'_{rs}+C'_{00}:r,s\in\{0,1\}\}.
\end{align*}
The corner equations have the separated form
\begin{equation}\label{eq:paired}
  \langle C_{ij},C'_{rs}\rangle=h_{ij}+h'_{rs},
\end{equation}
where
\[
  h_{ij}=w_1^{(i)}w_2^{(j)}+\omega_1^{(i)}+\omega_2^{(j)},
  \qquad
  h'_{rs}=w_3^{(r)}w_4^{(s)}+\omega_3^{(r)}+\omega_4^{(s)}.
\]
Taking a difference in each index pair gives
\[
  \langle C_{ij}+C_{00},C'_{rs}+C'_{00}\rangle=0.
\]
Thus $L\perp L'$, and nondegeneracy gives
\begin{equation}\label{eq:dim_sum}
  \dim_k L+\dim_k L'\le4.
\end{equation}
Summing \eqref{eq:paired} over $i,j$ at $r=s=0$ also gives
\begin{equation}\label{eq:face}
  \left\langle\sum_{i,j}C_{ij},C'_{00}\right\rangle
  =(w_1^{(0)}+w_1^{(1)})(w_2^{(0)}+w_2^{(1)}).
\end{equation}

We claim that $\dim_k L\ge2$, and that $e\in L$ whenever
$\dim_k L=2$. There are three cases.

If neither $u_1$ nor $u_2$ varies, then \eqref{eq:distinct} forces
both $w_1$ and $w_2$ to vary. But $\sum_{i,j}C_{ij}=0$, contradicting
\eqref{eq:face}.

If exactly one of $u_1,u_2$ varies, the difference in the other
coordinate is a nonzero multiple of~$e$. The difference in the varying
coordinate has a nonzero $K$-component. Since
$\sum_{i,j}C_{ij}=0$, these two differences span~$L$. Hence
$\dim_k L=2$ and $e\in L$.

Suppose that both $u_1,u_2$ vary, and let
$A_{ij}=Q^{-1}(u_1^{(i)}+u_2^{(j)})$ be the $K$-components of $C_{ij}$.
If the four points $A_{ij}$ are affinely independent, then
$\dim_k L=3$. Otherwise, Lemma~\ref{lem:rectangle} gives
\[
  \sum_{i,j}A_{ij}=0,
  \qquad \sigma=\sum_{i,j}\Tr(u_1^{(i)}A_{ij})\ne0.
\]
Consequently, $\sum_{i,j}C_{ij}=\sigma e\in L$. The projection of~$L$
onto~$K$ is one-dimensional and nonzero, while the kernel of that
projection is $ke$. Thus $\dim_k L=2$ and $e\in L$.

The same claim holds for $L'$. By \eqref{eq:dim_sum}, both spaces have
dimension~$2$ and contain~$e$. Their orthogonality gives the
contradiction
\[
  0=\langle e,e\rangle=1.\qedhere
\]
\end{proof}

\section{Consequences}\label{sec:consequences}

We first pass from powers of~$q^5$ to arbitrary side lengths.

\begin{corollary}\label{cor:lower}
For every integer $n\ge1$,
\[
  f_4(n)\ge2^{-19/20}n^{19/5}.
\]
Moreover, $\displaystyle\limsup_{n\to\infty}f_4(n)/n^{19/5}\ge1$.
\end{corollary}

\begin{proof}
If a box-free set $S_0\subseteq[N]^4$ is restricted to independently
chosen, uniformly random $n$-element subsets of the four parts,
where $n\le N$, its expected size is $|S_0|(n/N)^4$. Hence
\[
  f_4(n)\ge |S_0|(n/N)^4.
\]
For $n\ge32$, let $q$ be the largest power of~$2$ with $q^5\le n$,
and write $n=\mu q^5$, where $1\le\mu<32$. Embedding the
construction at~$q$ and sampling the construction at~$2q$ give
\[
  f_4(n)\ge \max\left\{\mu^{-19/5},\frac{\mu^{1/5}}2\right\}
                  n^{19/5}\ge2^{-19/20}n^{19/5},
\]
because the two coefficients agree at $\mu=2^{1/4}$, with the first
decreasing and the second increasing.
For $2\le n<32$, sampling at $q=2$ gives
$f_4(n)\ge n^4/2>2^{-19/20}n^{19/5}$, while $f_4(1)=1$.
Finally, $f_4(q^5)\ge q^{19}$ gives the limsup assertion.
\end{proof}

For $c\ge1$, let $R_4(c)$ be the least integer~$n$ such that every
$c$-coloring of $[n]^4$ contains a monochromatic box~\cite{CFP11}.

\begin{corollary}\label{cor:coloring}
For every power $q\ge2$ of~$2$, we have $R_4(q)>q^5$.
For every integer $c\ge2$,
\[
  R_4(c)>\frac{c^5}{32}.
\]
In particular, $R_4(2)\ge33$.
\end{corollary}

\begin{proof}
The coloring in Theorem~\ref{thm:main} gives $R_4(q)>q^5$.
For $c\ge2$, let $q$ be the largest power of~$2$ not exceeding~$c$.
Since $q>c/2$, monotonicity gives
$R_4(c)\ge R_4(q)>q^5>c^5/32$.
\end{proof}

The known exponents $19/5$ and $31/8$ differ by $3/40$.
It remains open whether the upper exponent is attainable, whether
the construction can be extended to more dimensions, and what the
correct growth rate of $R_4(c)$ is.

\end{document}